\documentclass[aos,MSNbibl,nameyear,dvips]{arximspdf}
\usepackage{graphicx}
\doi{10.1214/13-AOS1175B} 
\referstodoi{10.1214/13-AOS1175}
\volume{42}
\issue{2}
\pubyear{2014}
\firstpage{478}
\lastpage{482}

\makeatletter
\newcommand{\AR}{\mathrm{AR}}
\def\supp{\operatorname{supp}}
\makeatother

\begin{document}
\begin{frontmatter}

\title{Discussion: ``A significance test for the lasso''}
\runtitle{Discussion}

\begin{aug}
\author[a]{\fnms{T. Tony} \snm{Cai}\corref{}\thanksref{t1}\ead[label=e1]{tcai@wharton.upenn.edu}}
\and
\author[b]{\fnms{Ming} \snm{Yuan}\thanksref{t2}\ead[label=e2]{myuan@wisc.edu}}
\runauthor{T. T. Cai and M. Yuan}
\affiliation{University of Pennsylvania and University of Wisconsin--Madison}
\address[a]{Department of Statistics\\
The Wharton School\\
University of Pennsylvania\\
Philadelphia, Pennsylvania 19104\\
USA\\
\printead{e1}} 
\address[b]{Department of Statistics\\
University of Wisconsin--Madison\\
Madison, Wisconsin 53706\\
USA\\
\printead{e2}}
\pdftitle{Discussion of ``A significance test for the lasso''}
\end{aug}
\thankstext{t1}{Supported in part by NSF Grant DMS-12-08982 and NIH
Grant R01 CA127334.}
\thankstext{t2}{Supported in part by NSF Career Award DMS-08-46234 and
FRG Grant DMS-12-65202.}

\received{\smonth{12} \syear{2013}}



\end{frontmatter}

We congratulate the authors for an interesting article and an
innovative proposal to testing the significance of the predictor
variables selected by the Lasso. There is much material for thought and
exploration. Research on high-dimensional regression has been very
active in recent years, but most of the efforts have so far focused on
estimation. Despite the popularity of the Lasso as a variable selection
technique, the problem of making valid inference for a model chosen by
the Lasso is largely unsettled.
The current paper pinpoints some of the challenges in making valid
inference in the high-dimensional setting and presents a
thought-provoking approach to address them.

Following the notation used in the paper, let $A$ be the model selected
at the $k$th step of either the Lasso or forward stepwise regression
and $j$ be the index of the variable to be added in the next step. This
paper considers the problem of testing the null hypothesis that the
underlying model corresponding to the true regression coefficient
vector $\beta^\ast$ is nested in the current selected model, that is,
\[
H_0\dvtx  \supp\bigl(\beta^\ast\bigr)\subseteq A.
\]

As pointed out in the paper, a classical approach to testing two fixed
nested models $A$ and $A\cup\{j\}$ is the chi-squared test, which is
based on the test statistic
\[
R_j=(\mathrm{RSS}_A-\mathrm{RSS}_{A\cup\{j\}})/
\sigma^2
\]
and compares it to the quantile of the $\chi^2_1$ distribution. The
test fails, as noted, when applying to the forward stepwise regression
or the Lasso in a vanilla fashion because it fails to account for the
fact that neither $A$ nor $\{j\}$ is fixed. The randomness of $A$ can
be addressed using a conditional argument as suggested by the authors.
The effect of the way that the new index $j$ is selected is more
subtle. The seemingly lack of a remedy to this problem motives the
authors to focus on the Lasso and to propose the so-called covariance
test statistic
%
%
\begin{eqnarray} \label{Tkdecomp}
T_k&=& \bigl( \bigl\langle y, X\hat{\beta}(\lambda_{k+1})
\bigr\rangle- \bigl\langle y, X_A\tilde{\beta}_A(
\lambda_{k+1}) \bigr\rangle\bigr)/\sigma^2
\nonumber\\[-8pt]\\[-8pt]
&=&R_j-\lambda_{k+1} \bigl( \bigl\langle s_{A\cup\{j\}},
\hat{\beta}^{\mathrm{LS}}_{A\cup\{j\}} \bigr\rangle- \bigl\langle
s_{A},\hat{\beta}^{\mathrm{LS}}_{A} \bigr\rangle\bigr)/
\sigma^2,\nonumber
\end{eqnarray}
where $s_A$ and $s_{A\cup\{j\}}$ are, respectively, the vector of signs
of the nonzero regression coefficients for the Lasso at the $k$th and
$(k+1)$st steps, and $\hat{\beta}^{\mathrm{LS}}_{M}=(X_M^\top
X_M)^{-1}X_M^\top y$ is the least squares estimate under model $M$. In
effect, the second term on the right-hand side of (\ref{Tkdecomp})
can be viewed as a correction factor to account for the fact that the
next index $j$ is not fixed, but selected through the penalized $\ell
_1$ minimization. It is shown in the present paper that under $H_0$,
the limiting null distribution of $T_k$ is either $\operatorname{Exp}(1)$ or
stochastically smaller than $\operatorname{Exp}(1)$, and the paper proposed a test for
the null hypothesis $H_0$ based on this fact.

In this discussion, we introduce and explore a perhaps simpler and more
generic correction factor whose simplicity makes it an appealing
alternative to the current proposal. Furthermore, it can be easily
extended to other settings such as logistic regression and Cox
proportional hazards regression.

\subsection*{An alternative test}

Our proposal is based on the observation that for a given subset $A$,
the next selected index $j$ is not an arbitrary index in $A^c$. It is
instructive to first look at the case of orthogonal design where it is
clear that for both forward stepwise regression and the Lasso, $j$ can
be identified with
\[
R_j=\max_{m\in A^c} R_m.
\]
As a result, although for a fixed index $m\in A^c$, $R_m$ is a $\chi
^2_1$ distributed random variable, $R_j$, which is the maximum of $R_m$
for all $m\in A^c$, is not $\chi^2_1$ distributed. Note that,
conditioning on the design matrix $X$, $R_m$'s are independent $\chi
^2_1$ random variables. Therefore, the conditional distribution of
$R_j$ given $X$ can be easily deduced from the distribution of the
maxima of independent Gaussian random variables [see, e.g., \citet{HF06}]. In particular, in a high-dimensional setting where $p$
is large and $|A|$ is relatively small, the null distribution of $R_j$
can be well approximated by a Gumbel distribution (of type I). More
specifically, it can be shown that
%
%
\begin{equation}
\label{eqasy} 
\qquad R_j-2\log\bigl(\bigl|A^c\bigr|
\bigr)+\log\log\bigl(\bigl|A^c\bigr| \bigr) \stackrel{d} {\to} \operatorname{Gumbel}(-
\log\pi, 2)\qquad\mbox{as } p\to\infty,
\end{equation}
where the distribution function of a random variable $G$
following\break
$\operatorname{Gumbel}(-\log\pi, 2)$ is given by
\[
\mathbb{P}(G\le x)=\exp\bigl(-\exp\bigl(-(x+\log\pi)/2 \bigr) \bigr).
\]
This motivates us to consider the following test statistic:
%
%
\begin{equation}
\label{eqnewtest} 
\widetilde{T}_k=R_j-2\log
\bigl(\bigl|A^c\bigr| \bigr)+\log\log\bigl(\bigl|A^c\bigr| \bigr)
\end{equation}
and compare $\widetilde{T}_k$ with the quantile of $\operatorname{Gumbel}(-\log\pi, 2)$
distribution for testing the null hypothesis $H_0$. More specifically,
for any given $0<\alpha< 1$, we will reject $H_0$ at the $\alpha$
level if and only if $\widetilde{T}_k\ge q_{1-\alpha}^G$ where
$q_{1-\alpha}^G$ is the $1-\alpha$ quantile of $\operatorname{Gumbel}(-\log\pi, 2)$.

\begin{figure}

\includegraphics{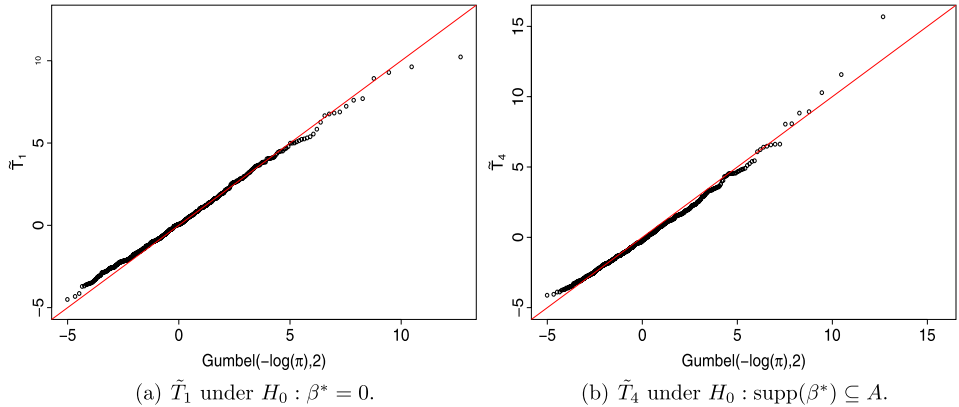}

\caption{Comparisons of the empirical distributions with the reference
distribution for $\widetilde{T}_k$ under the orthogonal design.}\label{figqqplot}
\end{figure}

To illustrate the accuracy of the reference distribution, we first
repeated the experiment considered in the paper with $n=100$
observations and $p=50$ variables under the orthogonal design. When the
true model is $\beta^\ast=0$ and, therefore, the null hypothesis
holds, we computed $\widetilde{T}_1$ for $500$ simulated datasets. The Q--Q
plot of the observed $\widetilde{T}_1$ versus its reference distribution
$\operatorname{Gumbel}(-\log\pi, 2)$ is given in the left panel of Figure~\ref
{figqqplot}. Similarly, the right panel of Figure~\ref{figqqplot}
gives the Q--Q plot for $\widetilde{T}_4$, again computed from 500
simulated datasets, when $\beta^\ast=(6,6,6,0,\ldots)^\top$.

%

\begin{figure}

\includegraphics{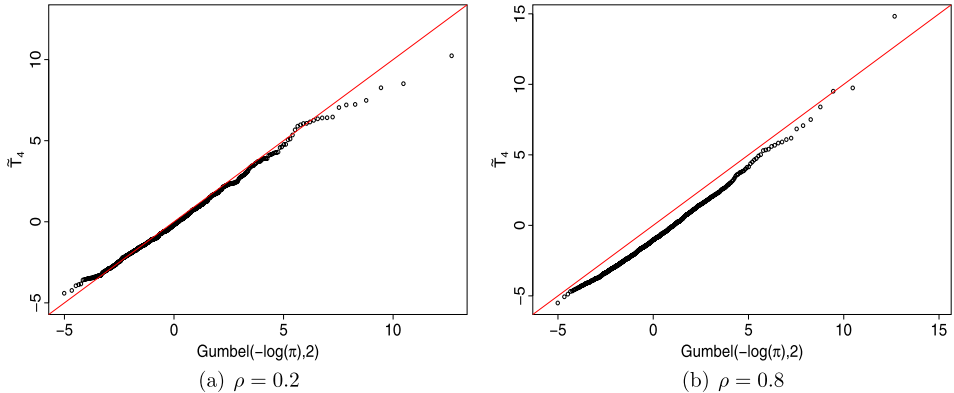}

\caption{Comparisons of the empirical distributions with the reference
distribution for $\widetilde{T}_4$ under the null $H_0\dvtx \operatorname{supp}(\beta
^\ast)\subseteq A$ with the $\AR(1)$ design.}
\label{figcorr}
\end{figure}

The strength of $\widetilde{T}_k$ comes from the robustness of its
limiting distribution under correlated designs. When $X^\top X\neq I$,
$R_m$'s are no longer independent but they are still marginally $\chi
^2_1$ distributed random variables. The distribution of $R_j=\max_{m\in
A^c}R_m$ again can be deduced from that of the maxima of a
Gaussian process. In particular, it can be shown that the limiting
Gumbel distribution given by (\ref{eqasy}) continues to hold under
fairly weak conditions on the dependence structure [see, e.g.,
\citet{LLR83}]. To verify the accuracy of
the Gumbel approximation under dependency, we repeated the previous
example with $\beta^\ast=(6,6,6,0,\ldots)^\top$. But instead of the
orthogonal design, the design matrix is now generated from a
multivariate normal distribution with mean zero and covariances $\operatorname{cov}(X_i, X_j)=\rho^{|i-j|}$. The left panel of Figure~\ref{figcorr}
corresponds to $\rho=0.2$ and right panel to $\rho=0.8$, both
suggesting that the limiting distribution $\operatorname{Gumbel}(-\log\pi, 2)$
continues to provide a reasonable approximation to the null
distribution of $\widetilde{T}_k$. In contrast, numerical results show
that the distribution of $T_k$ could deviate significantly from the
reference distribution $\operatorname{Exp}(1)$ under the correlated designs, and thus
comparing it to $\operatorname{Exp}(1)$ could be rather conservative in the correlated case.

%

\subsection*{General nonlinear $\ell_1$ regularization problems}

The advantages of the test statistic $\widetilde{T}_k$ proposed in (\ref
{eqnewtest}) are in its simplicity and generality. The correction
factor utilized by $\widetilde{T}_k$ depends only on the number of
remaining variables, and is\vadjust{\goodbreak} straightforward to evaluate. This makes it
particularly appealing when considering extensions to more general
nonlinear $\ell_1$ regularization problems where the exact tuning
parameter $\lambda_{k+1}$ for the next knot is typically not known in
closed form and often has to be approximated using an iterative
procedure. On the other hand, the validity of the Gumbel distribution
as the reference distribution under $H_0$ remains when $R_j$ is
replaced by the commonly used likelihood ratio test statistics.

\begin{figure}

\includegraphics{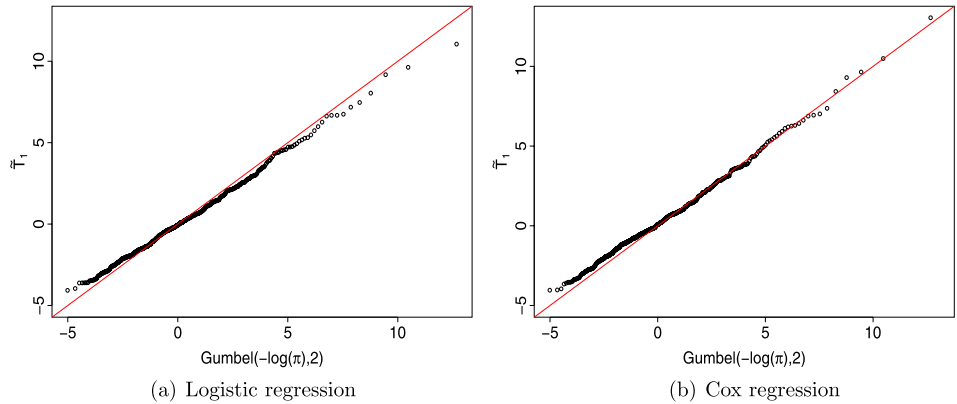}

\caption{Reference distribution for $\widetilde{T}_1$ under $H_0\dvtx \beta
^\ast=0$ for logistic regression and Cox's proportional hazards
model.}\label{figglm}
\end{figure}

To illustrate this point, we consider a logistic regression model where
the true regression parameter is $\beta^\ast=0$. With $n=100$
observations on a binary response and $p=50$ covariates independently
generated from the standard normal distribution. Same as before, the
experiment was repeated for 500 times; the Q--Q plot of the resulting
statistic $\widetilde{T}_1$ with respect to the $\operatorname{Gumbel}(-\log\pi, 2)$
distribution is given in the left panel of Figure~\ref{figglm}. The
right panel of Figure~\ref{figglm} shows the results from a similar
experiment for Cox proportional hazards regression where the response
was generated from $\operatorname{Exp}(1)$ with 10\% censoring. In both cases, the
reference $\operatorname{Gumbel}(-\log\pi, 2)$ distribution provides a good
approximation to the null distribution of the test statistic $\widetilde{T}_1$.

%

\subsection*{Summary}
The Lasso is a popular method for the high-dimensional linear
regression and it is important to make statistical inference for a
model chosen by the Lasso. The authors raise intriguing inferential
questions in the paper and propose a novel method to addressing them.
The work sheds new insight on high-dimensional model selection using
the Lasso and will definitely stimulate new ideas in the future. The
alternative test based on the test statistic $\widetilde{T}_k$ given in
(\ref{eqnewtest}) merits further investigation for linear
regression, logistic regression and Cox proportional hazards
regression, under the high-dimensional setting. We thank the authors
for their interesting work.




\printaddresses

\end{document}